# The notion of Infinity within the Zermelo system and its relation to the Axiom of Countable Choice[1].

**Chailos George**
Department of Computer Science-Division of Mathematics
University of Nicosia, Cyprus
e-mail: chailos.g@unic.ac.cy

## Abstract

In this article we consider alternative definitions-descriptions of a set being *Infinite* within the primitive Axiomatic System of Zermelo, $Z$. We prove that in this system the definitions of sets being *Dedekind Infinite, Cantor Infinite* and *Cardinal infinite* are equivalent each other. Additionally, we show that assuming the Axiom of Countable Choice, $AC_{\aleph_0}$, these definitions are also equivalent to the definition of a set being *Standard Infinite*, that is, *of not being finite*. Furthermore, we consider the relation of $AC_{\aleph_0}$ (and some of its special cases) with the statement $SD$ "*A set is Standard Infinite if and only if it is Dedekind Infinite*". Among other results we show that the system $Z + SD$ is 'strictly weaker' than $Z + AC_{\aleph_0}$.

## 1. Introduction

In this section we briefly present the aims of this article and we provide an overview of its main results. Furthermore, we provide the necessary definitions and statements that appear in the development of our arguments.

The concept of *finiteness*, defined as done via natural numbers, is categorical, thus not problematic. If, however, the concept of being finite is considered to be more fundamental than that of a number, as strongly advocated, e.g., by Frege and Dedekind (where natural

---

numbers are defined as the cardinals of finite sets) problems arise. First of all the concept of being *finite* loses its absoluteness. Secondly, how should *finiteness* be defined? Dual problems arise in relation to the concept of being *infinite*. There are several alternative different definitions-descriptions of a set being infinite; the oldest one is due to Dedekind (see [3]). In our presentation we consider *Dedekind Infinite* ( *DInf* ), *Cantor Infinite* (*CInf* ) and *Cardinal infinite* ($W_{\aleph_0}$) alternative definitions of a set being infinite. These definitions and the *Standard* one of being *Infinite* ( *SInf* ) are all equivalent to each other within the axiomatic system $ZFC$ which includes the *Axiom of Choice, AC* , in full strength (for details we refer to the classical book in Set Theory by Jech [8]). In this paper we show that this is also true even in the weaker system $Z + AC_{\aleph_0}$ that it does not include the *Axiom of Replacement*, the *Axiom of Foundation* and in which the *Axiom of Choice* (in full strength) is replaced by the *Axiom of Countable Choice*, $AC_{\aleph_0}$. Additionally we show that if we drop $AC_{\aleph_0}$, the alternative definitions *DInf* , *CInf* , $W_{\aleph_0}$, are all equivalent each other in Zermelo system $Z$ [2], but they are not equivalent to *SInf* . Furthermore, we investigate how the "strength" of $Z + AC_{\aleph_0}$ is affected if we substitute $AC_{\aleph_0}$ with the statement $SD$ : "*A set is Standard Infinite iff it is Dedekind Infinite, that is, if there is a 'one to one' correspondence 'onto' one of its proper subsets*"[3]. In this line of thought, we show that assuming the consistency of $Z$ and using results from Set Theory and Model Theory, the system $Z + SD$ is 'strictly weaker' than $Z + AC_{\aleph_0}$. Thus, in mathematical areas that avoid any version of the Axiom of Choice we could consider $Z + SD$ instead of $Z + AC_{\aleph_0}$.

Here it is worth mentioning that K. Gödel [4] and P. Cohen [1] showed the independence of the Axiom of Choice from the system $ZF$ and hence from $Z$ [4]. Particularly, assuming $ZF$ is consistent, K. Gödel showed that the negation of the Axiom of Choice is not a theorem of $ZF$ by constructing an inner model (the constructible universe) which satisfies $ZFC = ZF + AC$; thus, showing that $ZFC$ is consistent. Assuming $ZF$ is consistent, P. Cohen employed the technique of *forcing*, developed for this purpose, to show that the Axiom of Choice ( $AC$ ) itself is not a theorem of $ZF$ by constructing a much more complex model which satisfies $ZF \neg AC$ ($ZF$ with the negation of Axiom of Choice added as axiom) and thus showing that $ZF \neg AC$ is consistent. Together these results establish that the Axiom of Choice is logically independent from $ZF$. It is worth noting that the assumption that $ZF$ is consistent is harmless because adding another axiom to an already inconsistent system cannot make the situation worse. Because of independence, the decision whether to use (any version of) the Axiom of Choice (or its negation) in a proof cannot be made by appeal to other axioms of set theory. The decision must be made on other grounds, some of which are discussed in the sequel.

---

[2] See Section 2 for the Axioms of the original Zermelo system $Z$ .
[3] See Chailos [2] where a logical paradox is resolved within $Z + SD$ .
[4] The Axiom of Choice is not the only significant statement which is independent of $ZF$ . For example, the Generalized Continuum Hypothesis ( $GCH$ ) is not only independent of $ZF$ , but also independent of $ZFC$ . However, $ZF$ plus $GCH$ implies $AC$ , making $GCH$ a strictly stronger claim than $AC$ , even though they are both independent of $ZF$ (for details see Jech[8]).

For the following recall that two sets $X, Y$ are called *equipotent* if and only if there exists a bijection from $X$ to $Y$ and we also write $card\, X = card\, Y$.

We give the formal 'standard' definition of a set to be finite.

**Definition 1.1 (Finite Set)**

*A set $X$ is finite if (and only if) there is a natural number $n \in \mathbb{N}$ such that there is a one to one correspondence between $X$ and $n$. Thus, $X$ is finite iff it is equipotent to some $n \in \mathbb{N}$.*

**Definitions 1.2**

(a) *We set $\omega$ the least inductive set guaranteed by the Axiom of Infinity (see section 2).*
(b) *A set $X$ is countable if and only if there is a bijection from $X$ to $\omega$, thus if $X$ is equipotent to $\omega$.*

Observe that according to the above definition $\omega$ is the least countable ordinal. The cardinality of $\omega$ and hence by 1.2 (b) of any countable set $X$ is denoted by $\aleph_0$. Thus, $card(X) = \aleph_0$.

**Definitions 1.3 (Various definitions of Infinite Sets)**

1. *SInf : A set $X$ is Standard Infinite iff it is not Finite* (as in Definition 1.1).
2. *$W_{\aleph_0}$ : A set $X$ is Cardinal Infinite iff it has cardinality at least $\aleph_0$, $card(X) \geq \aleph_0$.*
3. *CInf : A set $X$ is Cantor Infinite iff it contains a countable subset.*
4. *DInf : A set $X$ is Dedekind Infinite iff it contains a proper subset $B \subset X$ equipotent to it. Hence, $card(B) = card(X)$.*

The following statements shall be used extensively in the development of our main arguments.

**Statements 1.4**

1. $SD$ : If a set is Standard Infinite, then it is also Dedekind Infinite.
2. $SC$ : If a set is Standard Infinite then it is also Cantor Infinite.
3. $SW_{\aleph_0}$ : If a set is Standard Infinite then it is also Cardinal Infinit.

*Here we summarize our main results*

**(a)** In *Theorem 3.1* we show that in $Z + AC_{\aleph_0}$ the different definitions of Infinite sets, $SInf$, $DInf$, $CInf$, $W_{\aleph_0}$ are equivalent each other.

**(b)** *In Theorem 3.5*, which its proof is a variation of the proof of *Theorem 3.1*, we show that in Zermelo system $Z$ (without any version of the Axiom of Choice) the different definitions of Infinite sets $DInf$, $CInf$, $W_{\aleph_0}$ are equivalent each other. Henceforth, the Axiomatic Systems $Z + SD$, $Z + SC$ and $Z + SW_{\aleph_0}$ are equivalent each other. It is worth to note that in $Z$ the definitions $DInf$, $CInf$, $W_{\aleph_0}$ are not equivalent to $SInf$, since there is a model of $Z$ such that there exists a *Standard Infinite* set that is *not Dedekind Infinite*. Such a model is the Cohen's First Model A4 (see [6]).

**(c)** Lastly, we examine the relation of $Z + SD$ with $Z + AC_{\aleph_0}$. At first we show that the Axiom of Countable Choice on Finite Sets, $AC_{\aleph_0}^{fin}$, is a theorem of $Z + SD$. In our main result (Theorem 3.8) we show that $Z + SD$ (and hence, according to *Theorem 3.5,* any of its equivalent systems $Z + SC$, $Z + SW_{\aleph_0}$) is "strictly weaker" than $Z + AC_{\aleph_0}$. Specifically, assuming consistency of $Z$, we show that there is a model of $Z$ in which all $SD$, $SC$, $SW_{\aleph_0}$ hold, but $AC_{\aleph_0}$ fails.

From **(b), (c)**, as above, we conclude that since the equivalent definitions of Infinite sets $DInf$, $CInf$, $W_{\aleph_0}$ in $Z$ are more intuitively obvious than $AC_{\aleph_0}$, any proof that adopts any of these statements and avoids the use of $AC_{\aleph_0}$ is preferable. This is important in many areas of Mathematics, such as Foundation of Mathematics, Philosophy of Mathematics, as well as in Set Theory that avoids the Axiom of Choice. We note that the system $Z$ consists of constructive axioms[5] and it is indispensable from any set theoretical development that does not involve atoms. This is another reason that any of the systems $Z + SD$, $Z + SC$, $Z + SW_{\aleph_0}$ is preferable compared to $Z + AC_{\aleph_0}$. This is of special interest in discrete structures since in $Z$ without assuming any version of the Axiom of Choice the following is true.

**Lemma 1.5**: *Countable sets have the property of being Dedekind Infinite* (no use of $AC_{\aleph_0}$)

**Proof**: If $X$ is a countable set, by definition 1.2(b), $X = \{a_i\}_{i=0}^{\omega}$. Take now the set $Y = X \setminus \{a_0\}$ which is a proper subset of $X$ with cardinality $\aleph_0$. To see that $X$ is *Dedekind Infinite* consider the function $h: X \to Y \subset X$, where $h(a_i) = a_{i+1}$ for every $i \in \omega$. Clearly, $h$ is 'one to one', since if $a_i \neq a_j$, $h(a_i) = a_{i+1} \neq a_{j+1} = h(a_j)$. In addition, $h$ is also 'onto', since $Range(h) = Y$. Thus, $h: X \to Y \subset X$ is a bijection from $X$ onto a proper subset of it, and thus $X$ is a *Dedekind Infinite* set.

## 2. The Set Theoretical Framework

Here we present and develop the necessary set theoretical framework required for proving the main theorems of the paper that are presented in section 3.

**The Axioms of $Z + AC_{\aleph_0}$**

Suppose that $\varphi(x_1, x_2, \ldots x_n)$ is a $n-$place well-formed formula (wff) in the formal language of set theory $\mathcal{L}$, which is a first order language that in addition to the usual symbol of equality it involves a symbol $\in$ for a binary predicate called *membership*.

***1. Axiom of Extensionality***: *If two sets have the same number of elements, then they are identical.*

$\forall x \forall y \forall z [(z \in x \Leftrightarrow z \in y) \Rightarrow x = y]$

This axiom does not tell us how we construct sets, but it tells what sets are: *A set is determined by its elements* (and it is atomless).

The Axioms 2-7 are all constructive axioms

***2. Null Set Axiom***: *There is an empty set, one which contains no element.* $\exists x \forall y \neg (y \in x)$.

---

[5] Apart from the *Axiom of Extensionality* which tells us what sets are (see section 2).

*Using the Axiom of Extensionality one can easily see that the empty set is unique. We write $\emptyset$ for the empty set.*

**3. Pair Set Axiom:** *If $a$ and $b$ are sets, then there is a set $\{a\}$ whose only element is $a$ and there is a set $\{a,b\}$ whose only elements are $a$ and $b$.*

$$\forall x \forall y \exists z \forall w (w \in z \Leftrightarrow w = x \vee w = y)$$

**4. Union Set Axiom:** *If $a$ is a set, then there is a set $\bigcup a$, the union of all the elements of $a$, whose elements are all the elements of elements of $a$.*

$$\forall x \exists y \forall z [z \in y \Leftrightarrow \exists w (w \in x \wedge z \in w)]$$

**5. Axiom of Infinity:** *There is a set which has the empty set, $\emptyset$, as an element and which is such that if $a$ is an element of it, then $\bigcup \{a, \{a\}\}$ (or $a \cup \{a\}$) is also an element of it. Such a set is called an inductive set:*

$$\exists x [\emptyset \in x \wedge \exists y (y \in x \Rightarrow \exists z (z \in x \wedge \forall w (w \in z \Leftrightarrow w \in y \vee w = y)))]$$

Consequently, this axiom guarantees the existence of a set of the following form:

$$\{\emptyset, \{\emptyset\}, \{\emptyset, \{\emptyset\}\}, \{\emptyset, \{\emptyset\}, \{\emptyset, \{\emptyset\}\}\}, \ldots\}$$

**6. Power Set Axiom:** *For any set $x$ there exists a set $y = P(x)$, the set of all subsets of $x$.*
That is, $\forall x \exists y \forall u [u \in y \Leftrightarrow \forall z (z \in u \Rightarrow z \in x)]$.

**7. Axiom Schema of (restrictive) Comprehension:** *If $a$ is a set and $\varphi(x, u_1, u_2, \ldots u_n)$ is a wff in $\mathcal{L}$, where the variable $x$ is free and $u_1, u_2, \ldots u_n$ are parameters (in a model of set theory from the axioms stated up to now), then there exist a set $b$ whose elements are those elements of $a$ that satisfy $\varphi$.*

$$\forall x \exists y \forall z [z \in y \Leftrightarrow z \in x \wedge \varphi(z, u_1, u_2, \ldots u_n)]$$

What this essentially means is that for any set $a$ and any wff $H[x]$ with one free variable $x$ and with parameters in a model of set theory (by the axioms defined as of now), there exist a set, unique by the *Axiom of extensionality*, whose elements are precisely those elements of $a$ that satisfy $H$. We denote this set by $b = \{x \in a : H[x]\}$. Hence, this axiom yields to subsets of a given set. Therefore, as a consequence of it, if $f : a \to b$ is a function with domain $a$, and codomain $b$, then $Range(f)$ is a well-defined set that is a subset of $b$.

The above seven axioms constitute the system $Z$ (Zermelo).

8. **Axiom of Countable Choice**: $AC_{\aleph_0}$ : *Suppose $b$ is any set and $P \subseteq \omega \times b$ any binary relation between natural numbers and members of $b$, then:*

$$\forall b \left[ (\forall n \in \omega)(\exists y \in b) P(n, y) \Rightarrow (\exists f : \omega \to b)(\forall n \in \mathbb{N}) P(n, f(n)) \right].$$

Essentially, the above axiom states that every Countable Family of nonempty sets has a choice function. That is, if $\mathcal{F} = \{A_i\}_{i \in \omega}$ is a countable family of nonempty sets, then there is $f : \mathcal{F} \to \bigcup \mathcal{F}$ such that $f(A_i) \in A_i$ for every $A_i \in \mathcal{F}$. In contrast to the full Axiom of Choice ($AC$) that demands the existence of choice functions $f : \mathcal{F} \to \bigcup \mathcal{F}$ for arbitrary families of nonempty sets, the Axiom of Countable Choice, $AC_{\aleph_0}$, justifies only a sequence of independent choices from arbitrary countable families of nonempty sets (see Ch.8 in [10]). This is equivalent to the statement "If $\langle X_i \rangle_{i \in \omega}$ is a countable indexed set of nonempty sets, then there exists a *(choice) function* $f : \omega \to \bigcup_{i \in \omega} X_i$ such that $\forall i \in \omega$, $f(i) \in X_i$. That is $\prod_{i \in \omega} X_i \neq \emptyset$." (See Ch.8, Theorem 1.3 of [7]). In light of the above, $AC_{\aleph_0}$ is more intuitively natural and is in concord with the techniques and the general spirit of mathematics that are based on Discrete/Countable structures (see Lemma 1.5).

Now we present and in some cases we prove necessary results from set theory.

At first we introduce the *Hartogs* number of any set.

**Definition 2.1.** The *Hartogs* number of any set $A$ is <u>the least ordinal number</u> $\alpha = h(A)$ which is not equipotent to any subset of $A$. That is, there is no injection $f : \alpha \to A$.

For clarification and further details see Ch7.1. of [7]. There, it is shown that the *Hartogs* number exists for any set and it is an initial ordinal (and thus it is a cardinal).

**Remark 2.2:**. Note that if $\lambda < h(A)$ [6], then there is an injection $g : \lambda \to A$. Hence, if $A$ is well-orderable set, $h(A)$ is the least cardinal greater than $card(A)$ and thus, $\lambda \leq card(A)$ (where here we consider $card(A)$ as the ordinal with which the cardinality $card(A)$ is identified).

---

[6] where $<$ is the standard ordering on ordinals.

An important concept for linearly (or totally) ordered sets, needed for the sequel, is that of *cofinality*.

With $(A, \prec)$ we denote a linearly ordered set.

**Definiton 2.3:** Let $(A, \prec)$ be a linearly ordered set. The set $(B, \prec) \subseteq (A, \prec)$ is called a cofinal subset in $(A, \prec)$, if it is <u>well ordered</u> [7] and such that $\forall a \in A\ \exists b \in B$ with $a \preceq b$.

The order type of a *well ordered set* $(X, \prec)$, $type(X)$, is the unique ordinal number $(\alpha, <)$ which is isomorphic to $(X, \prec)$ [8].

**Definition 2.4:** The cofinality of $(A, \prec)$ is the smallest possible order type of such cofinal set $(B, \prec)$ in $(A, \prec)$. That is, $cf(A) = \min\{type(B) : B \subseteq A \text{ and } B \text{ cofinal in } A\}$.

The following theorem is central for the development of our results and shall be used extensively.

**Theorem 2.5:** (Central theorem for existence of cofinal sets ): In $ZFC$, for every linearly ordered set $(A, \prec)$ there exists a <u>transfinite increasing sequence</u> $\hat{B}$, where its range $(B, \prec)$ is cofinal in $(A, \prec)$. Furthermore, $type(B) \leq card(A)$. (Where here we consider $card(A)$ as the ordinal with which the cardinality $card(A)$ is identified).

**Proof of Theorem 2.5:** We construct a transfinite increasing sequence $\hat{B} = \langle a_\xi : \xi < \lambda \rangle$ for some $\lambda < h(A)$, where $type(B) \leq card(A)$.

Fix $b \notin A$ (such an element exists, since $A$ is a set). By the Axiom of Choice there exists a choice function $g$ for the power set of $A$, $P(A)$.

We define the sequence $\langle a_\alpha : \alpha < h(A) \rangle$ by transfinite recursion (see Ch. 6, Theorem 4.4 in [7]) as follows: Given $\langle a_\xi : \xi < \alpha \rangle$, we consider two cases:

---

[7] And hence it can be viewed as an increasing sequence.
[8] For a proof, using the *Axiom of Replacement*, that such an ordinal exists see Theorem 6.3.1 in [7].

<u>Case 1</u>: If $b \neq a_\xi$, for all $\xi < \alpha$ and $A_\alpha = \{a \in A : a_\xi \prec a \text{ holds for all } \xi < \alpha\} \neq \emptyset$, we let $a_\alpha = g(A_\alpha)$.

<u>Case 2</u>: Otherwise we let $b = a_\alpha$.

Observe that from the recursive construction of $\langle a_\xi : \xi < \alpha \rangle$ this sequence is increasing.

*Claim*: There exists $\alpha < h(A)$ such that $a_\alpha = b$.

*Proof of Claim*: If for all $\xi < h(A)$ the terms $a_\xi$ of the increasing sequence, as constructed above, are such as $a_\xi \neq b$, then $a_\xi \in A$. Hence $\varphi : h(A) \to A$ will be an injection. This contradicts the definition of *Hartogs* number (see definition 2.1). Thus, there exists $\alpha < h(A)$ such that $a_\alpha = b$. This proves the claim. □

Now let $S = \{\alpha < h(A) : a_\alpha = b\}$. By claim, $S \neq \emptyset$ and since $S$ is a well-ordered set, it has a <u>least element</u> $\lambda = \min S$. From this and the recursive construction of the sequence $\hat{B} = \langle a_\xi : \xi < \lambda \rangle$ we get that $A_\lambda = \emptyset$ and thus, $B = \text{Range}\, \hat{B} \subseteq A$. Hence, $(B, \prec) \subseteq (A, \prec)$ is a well-ordered set with $\lambda = type(B) \leq card(A)$. (**) (See remark 2.2.)

Moreover, since $\lambda < h(A)$ with $A_\lambda = \emptyset$, there is no $a \in A$ such that for all $\xi < \lambda$, $a_\xi \prec a$. Since $(A, \prec)$ is linearly ordered, $\forall a \in A \, \exists \xi < \lambda$ with $a_\xi \in B$ such that $a \preceq a_\xi$. Thus, $\hat{B} = \langle a_\xi : \xi < \lambda \rangle$ is an increasing sequence such that $B = \text{Range}\, \hat{B} \subseteq A$ is cofinal in $(A, \prec)$ with $\lambda = type(B) \leq card(A)$. This concludes the proof of the theorem □

It is worth to restate the upshot of the above proof:

**Corollary 2.6:** For every linearly ordered set $(A, \prec)$ there exists a <u>transfinite increasing sequence</u> $\hat{B} = \langle a_\xi : \xi < \lambda \rangle$ cofinal in $(A, \prec)$, such that $\lambda < h(A)$.

Note that $\hat{B} = \langle a_\xi : \xi < \lambda \rangle$, which is cofinal in $(A, \prec)$ can be viewed as an increasing function $f : (\lambda, <) \to (A, \prec)$ such that $f(\xi) = a_\xi$ for $\xi < \lambda$. In such a case we say that $f : (\lambda, <) \to (A, \prec)$ is a cofinal map from $(\lambda, <)$ to $(A, \prec)$. In particular, for ordinals $\alpha$, $\beta$ (and in general for well-ordered sets) the above analysis justifies the following definition.

**Definition 2.7:** (cofinal maps): Let $\alpha$, $\beta$ ordinals with $\beta \leq \alpha$. A map $f : \beta \to \alpha$ is cofinal in $\alpha$ if it is increasing and is such that $\forall \gamma < \alpha \; \exists \delta < \beta$ such that $\gamma \leq f(\delta)$.

According to all of the above (see 2.4, 2.5, 2.6, 2.7) we have the following well defined equivalent definition of cofinality of ordinals.

**Definition 2.8** (cofinality of ordinals): The cofinality of an ordinal $\alpha$ is defined to be the minimal ordinal $\beta$ such that there exists a cofinal map $g : \beta \to \alpha$. That is:

$$cf(\alpha) = \min\{\beta : \exists g : \beta \to \alpha \text{ which is cofinal in } \alpha\}.$$

Having these equivalent definitions of cofinality[9] in our disposal, we can easily obtain results for the cofinality on ordinals that are required in the sequel.

**Facts about the notion of cofinality**

(a) $cf(0) = 0$ (Immediate from definition 2.4).

(b) Since for any ordinal $\gamma$, $id : \gamma \to \gamma$ is a cofinal map in $\gamma$, we get that $cf(\gamma) \leq \gamma$.

(c) An ordinal $\alpha$ is a successor if and only if $\alpha = \beta \cup \{\beta\}$ and hence if and only if $cf(\alpha) = 1$. To see this consider the cofinal map $f : 1 \to \alpha$ such that $f(0) = \beta$.

(d) If $\alpha$ is a successor ordinal, then $cf(\aleph_\alpha) = \aleph_\alpha$. The proof is a consequence of the fundamental theorem of cardinal arithmetic. For details see Thm 9.2.4 in [7].

**Lemma 2.9**: Let $\alpha$ be a countable limit ordinal. Then $cf(\alpha) = \omega$

**Proof:** Since $\alpha$ is countable limit ordinal, there is an enumeration of $\alpha$. That is, $\alpha = \{\beta_n\}_{n \in \omega}$. Recursively define a map $f : \omega \to \alpha$ as: $f(0) = \min\{\beta_n\}_{n \in \omega}$ and $f(n) = \max\{f(n-1), \beta_n\}$. Then, by definition 2.7, it is easy to conclude that $f : \omega \to \alpha$ is a cofinal map in $\alpha$. Since $\omega$ is the least countable limit ordinal, by definition 2.8 $cf(\alpha) = cf(\omega) = \omega$. □

According to the above lemma we have, $cf(\omega + \omega) = cf(\omega \cdot \omega) = cf(\omega_\omega) = cf(\omega) = \omega$.

---

[9] For an alternative, but according to our above analysis, an equivalent approach to cofinality of ordinals we refer to Ch 9, Sec 2 in [7].

**Lemma 2.10:** If $\alpha, \beta, \gamma$ ordinals and $f: \beta \to \alpha$, $g: \gamma \to \beta$ are cofinal maps in $\alpha$ and $\beta$ respectively, then $f \circ g: \gamma \to \alpha$ is a cofinal map in $\alpha$.

**Proof**: Observe that since $f, g$ are increasing, $f \circ g$ is increasing. Now, since $f: \beta \to \alpha$ is cofinal, for every $\xi < \alpha$ there exists $\delta < \beta$ such that $\xi \leq f(\delta)$. Since $g: \gamma \to \beta$ is cofinal, there exists $\zeta < \gamma$ such that $\delta \leq g(\zeta)$. Because $f$ is an increasing map, $f(\delta) \leq f(g(\zeta))$. Putting everything together, $\xi \leq f(\delta) \leq f(g(\zeta))$ and thus, for every $\xi < \alpha$, there exists $\zeta < \gamma$, such that $\xi \leq f(g(\zeta))$. Therefore, $f \circ g: \gamma \to \alpha$ is a cofinal map in $\alpha$.

**Corollary 2.11**: If $f: \beta \to \alpha$ is a cofinal map in $\beta$, then $cf(\alpha) \leq cf(\beta)$.
(To see this, set $\gamma = cf(\beta)$ in the above lemma 2.10).

**Corollary 2.12:** For every ordinal $\alpha$, $cf(cf(\alpha)) = cf(\alpha)$.

**Proof:** In corollary 2.11 set $\beta = cf(\alpha)$. Then $cf(\alpha) \leq cf(cf(\alpha))$.

For the reverse inequality, set $\gamma = cf(\alpha)$ in $cf(\gamma) \leq \gamma$ to get $cf(cf(\alpha)) \leq cf(\alpha)$. □

**Lemma 2.13:** *If $\alpha$ is a limit ordinal, then $cf(\aleph_\alpha) = cf(\alpha)$*

**Proof**: Since $\alpha$ is a limit ordinal, by definition $\alpha = \sup\{\beta : \beta < \alpha\}$ and hence $\aleph_\alpha = \bigcup_{\beta < \alpha} \aleph_\beta$. The increasing map $f: \alpha \to \aleph_\alpha$ such that $f(\beta) = \aleph_\beta$[10], is cofinal in $\alpha$.

Now by corollary 2.11, $cf(\aleph_\alpha) \leq cf(\alpha)$. **(1)**

From (1) and fact (b), $cf(\aleph_\alpha) \leq \alpha$. Now let $f: \lambda \to \aleph_\alpha$ with $\lambda \leq \alpha$ be any cofinal map of $\lambda$ in $\aleph_\alpha$. Consider $g: \lambda \to \alpha$ such that $g(\xi) = \max\{\beta : \aleph_\beta \leq f(\xi)\}$.

*Claim*: $g: \lambda \to \alpha$ is cofinal in $\alpha$.

*Proof of Claim*: We first show that $g: \lambda \to \alpha$ is increasing. Indeed, let $\gamma < \delta < \lambda$, and let $\beta < \alpha$ such that $\aleph_\beta \leq f(\gamma)$. Since $f: \lambda \to \aleph_\alpha$ is increasing, $\aleph_\beta \leq f(\gamma) < f(\delta)$, and since $\aleph_\alpha = \bigcup_{\rho < \alpha} \aleph_\rho$, by definition of $g$ we conclude that $g(\gamma) < g(\delta)$. Thus, $g: \lambda \to \alpha$ is increasing. Furthermore, if $\gamma < \alpha$, then $\aleph_\gamma < \aleph_\alpha$ and since $f: \lambda \to \aleph_\alpha$ is cofinal in $\aleph_\alpha$, there is a $\zeta < \lambda$ such that $\aleph_\gamma \leq f(\zeta)$. Hence, by the maximality in the definition of $g : \lambda \to \alpha$ we conclude that $\gamma \leq g(\zeta)$. Therefore, for every $\gamma < \alpha$, there is $\zeta < \lambda$ such that $\gamma \leq g(\zeta)$. Hence, $g : \lambda \to \alpha$ is cofinal in $\alpha$. This proves the claim. □

---
[10] Note that for every $\beta < \gamma$ we have $\aleph_\beta < \aleph_\gamma$.

Note that by (1) we get that $cf(\aleph_\alpha) \leq cf(\alpha) \leq \alpha$. Since $f : \lambda \to \aleph_\alpha$ with $\lambda \leq \alpha$ is any cofinal map of $\lambda$ in $\aleph_\alpha$, by setting $\lambda = cf(\aleph_\alpha)$ in $g : \lambda \to \alpha$ we conclude by corollaries 2.11 and 2.12 that $cf(\alpha) \leq cf(cf(\aleph_\alpha)) = cf(\aleph_\alpha)$. **(2)**

From **(1), (2)**, $cf(\aleph_\alpha) = cf(\alpha)$ and the proof is now complete. □

The last definition of this section shall be used in the proof of *Theorem 3.8*.

**Definition 2.14**: An ordinal $\alpha$ is called *singular* iff $cf(\alpha) < \alpha$, otherwise it is called *regular*.

### 3. The Main Results

In this section we prove in detail the main results of this paper and we discuss their consequences.

In *Theorem 3.1* we show that in $Z + AC_{\aleph_0}$ the different definitions of Infinite sets, *SInf*, *DInf*, *CInf*, $W_{\aleph_0}$ are all equivalent each other.

*In Theorem 3.5*, which its proof is a variation of the proof of *Theorem 3.1*, we show that in $Z$ (without any version of the Axiom of Choice) a set is *Dedekind Infinite iff it is Cantor Infinite iff it is Cardinal Infinite*. Henceforth, the Axiomatic Systems $Z + SD$, $Z + SC$ and $Z + SW_{\aleph_0}$ are equivalent. It is worth to note that in $ZF$ (and hence in $Z$) the definitions *DInf*, *CInf*, $W_{\aleph_0}$ are not equivalent to *SInf*, since there is a model of $ZF$ such that there exists a *Standard Infinite* set that is *not Dedekind Infinite*. Thus, $SD$ (and hence $SC$, $SW_{\aleph_0}$) does not hold. Such a model is the Cohen's First Model A4 (see [6]).

We further study how $SD$ is related to $AC_{\aleph_0}$ and we show that $Z + AC_{\aleph_0} \vdash SD$[11]. That is, the statement $SD$ is a theorem of the system $Z + AC_{\aleph_0}$. Furthermore we show that $Z + SD \vdash AC_{\aleph_0}^{fin}$, where $AC_{\aleph_0}^{fin}$ is the Axiom of Countable Choice on Finite Sets.

---

[11] $X \vdash Y$ denotes that the statement $Y$ is a logical consequence of $X$. That is, there is a proof of $Y$ from $X$, and hence $Y$ is a theorem of $X$.

Lastly, we show in *Theorem 3.8*, using results from section 2 and model theory, that the system $Z + SD$ (and hence any of its equivalent-according to *Theorem 3.5*- systems, $Z + SC$, $Z + SW_{\aleph_0}$) is <u>strictly weaker</u> than $Z + AC_{\aleph_0}$.

In the following theorem we show that in $Z + AC_{\aleph_0}$ the different definitions of Infinite sets, $SInf$, $DInf$, $CInf$, $W_{\aleph_0}$ are all equivalent each other. In particular $Z + AC_{\aleph_0} \vdash SD$.

**Theorem 3.1:** *In $Z + AC_{\aleph_0}$ a set is Standard Infinite if and only if $card(X) \geq \aleph_0$, if and only if it is Cantor Infinite, if and only if it is Dedekind Infinite.* Thus, $SInf$, $DInf$, $CInf$, $W_{\aleph_0}$ are all equivalent each other.

**Proof:** For this, we show in (a)-(c) that $SInf \Rightarrow CInf \Rightarrow DInf \Rightarrow SInf$ and in (d) that $CInf \Leftrightarrow W_{\aleph_0}$.

(a): We show, $SInf \Rightarrow CInf$ : (That is, every standard infinite set contains a countable subset). The proof is rather simple if we assume the *Axiom of Replacement* and *the Axiom of Choice* in its full strength (see Appendix Lemma A1). <u>Here we prove it only from $Z + AC_{\aleph_0}$</u>.

Let $X$ be a *Standard Infinite* set. For every $n \in \mathbb{N}$, let $A_n$ be the set of all subsets of $X$ with cardinality $2^n$. Thus, the members of $A_n$ are sets of cardinality $2^n$. Since $X$ is a standard infinite set, by definitions 1.1 and 1.3, $\forall n \in \mathbb{N}, A_n \neq \emptyset$. Apply $AC_{\aleph_0}$ to the family $\mathcal{F} = \{A_i\}_{i \in \mathbb{N}}$ in order to obtain a choice function $\phi : \mathcal{F} \to \bigcup \mathcal{F}$, and via it a subsequence $\{B_n\}_{n \in \mathbb{N}} = \phi[\mathcal{F}]$, where each $B_n$ is a subset of $X$ of cardinality $2^n$. (It is worth to note that the *Axiom of Replacement* is not required since $\{B_n\}_{n \in \mathbb{N}} = \phi[\mathcal{F}]$ is the image in a predefined set of a well-defined function on a set[12].)

Observe that the sets $\{B_n\}_{n \in \mathbb{N}}$ are also <u>distinct</u>. We define recursively the sequence $\{C_n\}_{n \in \mathbb{N}}$ of <u>mutually disjoint</u> sets as follows (this is a standard technique in analysis):

$$C_0 = B_0, \ C_1 = B_1 \setminus C_0, \ \ldots, C_n = B_n \setminus \bigcup_{i=1}^{n-1} C_i, \ldots$$

Obviously, $1 \leq card(C_n) \leq 2^n$, since $\sum_{i=0}^{n-1} 2^i = 2^n - 1$. Now it is trivial to observe that the family $\mathcal{S} = \{C_n\}_{n \in \mathbb{N}}$ is a family of nonempty and mutually disjoint sets.

---

[12] See also section 2 the discussion following the *Axiom of Replacement* 7.

A second application of $AC_{\aleph_0}$ to the nonempty family $\mathcal{S}' = \{C_n\}_{n \in \mathbb{N}}$ guarantees the existence of a choice function $\psi : \mathcal{S}' \to \bigcup \mathcal{S}'$, and through this we obtain a set $C = \{c_n\}_{n \in \mathbb{N}} = \psi[\mathcal{S}']$, such that $c_i \neq c_j$, $\forall i \neq j \in \mathbb{N}$ (since $\{C_n\}_{n \in \mathbb{N}}$ is a sequence of mutually disjoint sets). This set $C = \{c_n\}_{n \in \mathbb{N}}$ is the countable set that establishes $X$ to be *Cantor Infinite*. (As done earlier in the proof, the *Axiom of Replacement* is not required since $C = \{c_n\}_{n \in \mathbb{N}} = \psi[\mathcal{S}']$) □

(b): We show, $CInf \Rightarrow DInf$ (that is, if a set is *Cantor Infinite* then it is *Dedekind Infinite*). Since $X$ is *Cantor Infinite* it contains a countable subset $B \subseteq X$. That is, there exists a bijection $f : \mathbb{N} \mapsto B$. Define the function $h : X \mapsto X$ as follows:

$h(x) = f(n+1)$ if $x = f(n)$, $x \in B$     and     $h(x) = x$, if $x \in X \setminus B$.

(Thus, the action of $h$ on $B$ is identified with the action of $f$ on the successor of each element of $\mathbb{N}$ that leaves the remaining elements of $X$ invariant.)

*Claim 1:* The function $h : X \mapsto X$ is 'one to one'.

*Proof of Claim 1:* Let $x, y \in X$ with $h(x) = h(y)$. In order for such a thing to hold we can easily see that either both $x, y \in B$, or $x, y \in X \setminus B$.

In the first case, where $x, y \in B$, we have $x = f(n)$, $y = f(m)$ for some $n, m \in \mathbb{N}$. From the definition of $h$, $h(x) = h(y) \Rightarrow f(n+1) = f(m+1)$, and since $f$ is 'one to one', $n+1 = m+1$ and hence $n = m$. Thus, $x = y$.

In the second case, where $x, y \in X \setminus B$, we conclude that $h(x) = h(y) \Rightarrow x = y$.

This proves claim 1. □

Let $D = h(X) \subseteq X$. Then the function $h : X \mapsto D$ from *claim 1* is 'one to one' and since $D = Range\ h$, then $h : X \mapsto D$ is also onto and hence a bijection. Thus, $X$ is equipotent to its subset $D = h(X)$. To complete the proof it is enough to show the following.

*Claim 2:* The set $D = h(X)$ is a proper subset of $X$.

*Proof of Claim 2:* Let $x = f(0)$. We claim that $x \notin D$. Indeed, if we suppose that $x \in D$, then $x = h(z)$ for some $z \in X$, and we have the following two cases.

(i) If $z \in B$, then $z = f(n)$ for some $n \in \mathbb{N}$. Hence by the definition of $h$, $x = f(0) = h(z) = f(n+1)$. Thus, $f(0) = f(n+1)$. Since $f$ is 'one to one', $0 = n+1$. This leads to a contradiction since $0 \in \mathbb{N}$ is not a successor.

(ii) If $z \in X \setminus B$, then $x = f(0) = h(z) = z$ and thus (from the definition of $f: \mathbb{N} \mapsto B$) $z \in B$. This leads again to contradiction (since by hypothesis $z \in X \setminus B$.)

Putting everything together, we conclude that there exists $x \in X \setminus h(X)$ and therefore $D = h(X)$ it is a proper subset of $X$. This proves claim 2 □

The proof of part (b) is now complete. □

**Remark 3.2**: *Observe that in (b) of the above proof we have shown that $CInf \Rightarrow DInf$ without using the Axiom of Countable Choice.*

(c): We show, $DInf \Rightarrow SInf$ (that is, if a set is *Dedekind Infinite*, then it is *Standard Infinite*- equivalently not finite). Suppose, for contradiction, that $X$ is a finite, yet a Dedekind infinite set. It is well known that in $Z$ (and hence in $Z + AC_{\aleph_0}$) the *Pigeonhole Principle* holds: "If a set is finite, then it does not contain a proper subset equipotent to it." (see Appendix Lemma A2). Hence, $X$ does not contain a proper subset equipotent to it. This contradicts the hypothesis of $X$ being *Dedekind Infinite* (and thus containing a proper subset equipotent to it). Therefore, $DInf$ implies $SInf$. □

**Remark 3.3**: *Observe that in (c) we have shown that $DInf \Rightarrow SInf$ without using the Axiom of Countable Choice.*

(d): We show, $CInf \Leftrightarrow W_{\aleph_0}$ (that is, a set is Cantor Infinite iff it is Carinal Infinite).

(i) $CInf \Rightarrow W_{\aleph_0}$. Suppose that $X$ is a *Cantor Infinite* set. Then $X$ contains a countable subset $A$ with $card(A) = \aleph_0$, and hence there is a bijection $\varphi: \mathbb{N} \to A \subseteq X$. Thus, $\aleph_0 = card(A) \leq card(X)$ (see also section 1). □

(ii) $W_{\aleph_0} \Rightarrow CInf$. Suppose that $X$ is a set such that $card(X) \geq \aleph_0$. Hence, there is an injection $g: \mathbb{N} \to X$. Thus, $g[\mathbb{N}] \subseteq X$ is a countable set and hence $X$ is Cantor Infinite. □

**Remark 3.4:** *Observe that in (d) we have shown $CInf \Leftrightarrow W_{\aleph_0}$ without using the Axiom of Countable Choice.*

Putting (a), (b), (c), (d) together we conclude that in $Z + AC_{\aleph_0}$,

$SInf \Leftrightarrow W_{\aleph_0} \Leftrightarrow CInf \Leftrightarrow DInf$. The proof of the theorem is now complete. □

**Theorem 3.5:** *In the axiomatic system Zermelo, $Z$, a set is Cardinal Infinite if and only if it is Cantor Infinite, if and only if it is Dedekind Infinite. Hence, $W_{\aleph_0}, CInf, DInf$ are all equivalent each other in $Z$.*

**Proof:** From remark 3.4, $CInf \Leftrightarrow W_{\aleph_0}$ in $Z$, and from remark 3.2, $CInf \Rightarrow DInf$ in $Z$. It remains to show that in $Z$, $DInf \Rightarrow CInf$. That is, if $X$ *is Dedekind Infinite*, then $X$ it is also *Cantor Infinite*. To this end, since $X$ is Dedekind infinite, it contains a proper subset $B$ equipotent to it. Thus, $\exists f : X \to B \subset X$ which is 'one to one' and 'onto' $B \subset X$. Consider $a \in X \setminus B$. Then, $a \notin Range(f)$. To show that $X$ is *Cantor Infinite* we construct recursively a well-defined countable (infinite) subset $Y$ of $X$ as follows:

Consider $f : X \to X$ as above. Then by (simple) recursion theorem (see 5.6 in [10]) there is a unique function $g : \mathbb{N} \to X$, such that:

$g(0) = a$ and $g(n+1) = f((g(n))$ for all $n \in \mathbb{N}$.

The set $Y = Range(g) = \{g(n)\}_{n \in \mathbb{N}}$ is a well-defined countable subset $Y \subseteq X$. Indeed, by construction, for any $n \in \mathbb{N}$, $g(n) = f^{(n)}(a)$ where $f^{(n)}$ denotes the $n^{th}$ recursive iteration of the action of $f$ on $a$. If (w.l.o.g) $\exists i < j$ such that $g(i) = g(j)$, then $f^{(i)}(a) = f^{(j)}(a)$ and hence $f^{(j-i)}(a) = a$ where $j - i \neq 0$. Thus, $a \in Range(f)$, and this contradicts the choice of $a$. Hence, $g$ is a 'one to one' well-defined function. This set $Y$ is the Infinite Countable set we are looking for. □

In the next theorem we prove that in $Z + SD$ the *Axiom of Countable Choice* on *Finite Sets*, $AC_{\aleph_0}^{fin}$, holds. Thus, $AC_{\aleph_0}^{fin}$ is a logical consequence (theorem) of $Z + SD$.

The main idea for the proof of the next theorem is from Theorem 2.12 in [5].

**Theorem 3.6** $Z + SD \vdash AC_{\aleph_0}^{fin}$.

**Proof:** Let $\mathcal{F} = \{X_n\}_{n \in \mathbb{N}}$ be a countable family of nonempty finite sets. We show that $\prod_{n \in \mathbb{N}} X_n \neq \varnothing$ and hence $AC_{\aleph_0}^{fin}$ holds[13]. For this, define $X = \bigcup_{n \in \mathbb{N}} (X_n \times \{n\})$. Then $X$ is not a finite set, equivalently it is *Standard Infinite* and thus by *SD* it is *Dedekind Infinite*. By Theorem 3.5 it is also *Cantor Infinite* and hence there is an injection $f : \mathbb{N} \mapsto X$. Since $\forall n \in \mathbb{N}$, $X_n$ is a finite set, using proof by contradiction, we can easily show that $M = \{n \in \mathbb{N} : f[\mathbb{N}] \cap (X_n \times \{n\}) \neq \varnothing\}$ is a standard infinite set.

*Claim 1:* $\prod_{n \in M} X_n \neq \varnothing$

*Proof of Claim 1:* Since $f[\mathbb{N}] \cap (X_m \times \{m\}) \neq \varnothing$, for every $m \in M$ define $n(m) = \min\{n \in \mathbb{N} : f(n) \in X_m \times \{m\}\}$. Thus, for every $m \in M$ there exists a unique $x_m \in X_m$ such that $f(n(m)) = (x_m, m)$. Henceforth, there exists an element $(x_m)_{m \in M} \in \prod_{m \in M} X_m$. Thus, $\prod_{n \in M} X_n \neq \varnothing$. □

According to the above, in order to complete the proof it is enough to show the following:

*Claim 2:* If $\mathcal{F} = \{X_n\}_{n \in \mathbb{N}}$ is a countable collection of nonempty finite sets such that it exists a standard infinite set $M \subseteq \mathbb{N}$ with the property $\prod_{m \in M} X_m \neq \varnothing$, then $\prod_{n \in \mathbb{N}} X_n \neq \varnothing$. (Hence, the *Axiom of Countable Choice* for *Finite Sets* holds.)

*Proof of Claim 2:* Consider the sequence of nonempty finite sets $\langle X_n \rangle_{n \in \mathbb{Z}_+}$ and define $Y_k = \prod_{n \leq k} X_n$. Inductively observe that for every $k \in \mathbb{Z}_+$ the sets $Y_k$ are nonempty and are finite. Now consider the sequence of nonempty finite sets $\langle Y_k \rangle_{k \in \mathbb{Z}_+}$. By hypothesis there exists a finite set $M \subseteq \mathbb{Z}_+$ with the property $\prod_{m \in M} Y_m \neq \varnothing$, and hence, there exists an element $(y_m)_{m \in M} \in \prod_{m \in M} Y_m$ with $y_m = (x_1^m, x_2^m, \ldots, x_m^m) \in \prod_{n \leq m} X_n$. Equivalently, $\forall \, 1 \leq i \leq m \in M$, $x_i^m \in X_i$. (*)

---

[13] See the discussion in section 2 that follows Axiom 8, *Axiom of Countable Choice*. There we refer explicitly to equivalent to $AC_{\aleph_0}$ statements.

Let $n \in \mathbb{N}$ and consider the set $Z_n = \{m \in M : n \leq m\}$. Since $M \subseteq \mathbb{Z}_+$ is a *Standard Infinite* set, by definition 1.1 we have that $\forall n \in \mathbb{Z}_+ \ \exists m \in M : n \leq m$. Thus, $\forall n \in \mathbb{Z}_+$ the set $Z_n$ is a nonempty subset of $\mathbb{N}$ and hence it has a minimum element. Now define $m(n) = \min Z_n = \min\{m \in M : n \leq m\}$. Since, $y_{m(n)} = \left(x_1^{m(n)}, x_2^{m(n)}, \ldots, x_{m(n)}^{m(n)}\right) \in \prod_{n \leq m(n)} X_n$ and $n \leq m(n)$, according to (*), $x_n^{m(n)} \in X_n$. The sequence $f : \mathbb{Z}_+ \mapsto \bigcup_{i \in \mathbb{Z}_+} X_i$, where $f(n) = x_n^{m(n)} \in X_n$, is a choice function for the family $\mathcal{F} = \{X_n\}_{n \in \mathbb{Z}_+}$, with $\left(x_n^{m(n)}\right)_{n \in \mathbb{Z}_+} \in \prod_{n \in \mathbb{Z}_+} X_n$ (see footnote 13). Therefore $\prod_{n \in \mathbb{Z}_+} X_n \neq \varnothing$, and since $X_0 \neq \varnothing$, then $\prod_{n \in \mathbb{N}} X_n \neq \varnothing$. This finishes the proof of claim 2

The proof of theorem is now complete. □

In our last result we show, using results from section 2 and model theory, that the system $Z + SD$ <u>is strictly weaker</u> than $Z + AC_{\aleph_0}$, as one could naturally conjecture from the preceeding results (see Theorems 3.1, 3.5 and 3.6). For this we show that there exists a model $\mathcal{M}$ of $Z$ where $SD$ holds, but $AC_{\aleph_0}$ does not hold. <u>Thus, the statement</u> $SD$, "Every *Standard Infinite* set is *Dedekind Infinite*", <u>is even weaker than the</u> *Axiom of Countable Choice,* $AC_{\aleph_0}$. Therefore, it could be naturally added to the axioms $Z$ (instead of $AC_{\aleph_0}$) in order to address problems in a framework that does not assumes any version of the Axiom of Choice. Henceforth, in this sense, we could develop a rather powerful axiomatic system in areas that any version of Axiom of Choice is unnatural[14]. This is of great interest in the area of constructive mathematics where the adoption of any version of Axiom of Choice is doubtful or unacceptable (see Richman[11]).

Recall, as stated earlier in the paper, that the statement $SD$ cannot be proved from $Z$ or $ZF$ alone, since there exists a model of $ZF$ (and hence of $Z$) in which there exists a standard Infinite set that is not Dedekind Infinite. Such a model is the Cohen's First Model A4 (see [6]).

**Statement** $\hat{W}_\lambda$: For every set $X$ and a cardinal $\lambda$, either $card(X) < \lambda$ or $card(X) \geq \lambda$.

---

[14] See Chailos [2] for resolution of Zeno paradoxes where $SD$ is used without any use of $AC_{\aleph_0}$.

Observe that $\hat{W}_{\aleph_0}$ is the statement that every set is either *Cardinal Finite* or *Cardinal Infinite* (see section 1). The following theorem is Theorem.8.6 of [9].

**Theorem 3.7:** *Let $\aleph_\alpha$ be a singular aleph. Then there exists a model of $ZF$ in which*

(i) *For each $\lambda < \aleph_\alpha$, $\hat{W}_\lambda$ holds.*

(ii) *$\hat{W}_{\aleph_\alpha}$ fails and $AC_{cf(\aleph_\alpha)}$ fails.*

**Theorem 3.8:** *In $Z$ the statement $SD$ is strictly weaker than $AC_{\aleph_0}$.*

**Proof**: From the proof of Theorem 3.1, $Z + AC_{\aleph_0} \vdash SD$. Thus, by soundness theorem any model of $Z$ that satisfies $AC_{\aleph_0}$ it also satisfies $SD$. Hence, in order to prove that the system $Z + SD$ is weaker than $Z + AC_{\aleph_0}$ it is enough to show that assuming the consistency of $Z$, there exists a model $\mathcal{M}$ of $Z$ that satisfies $SD$ in which $AC_{\aleph_0}$ fails. Since $\omega$ is a limit ordinal, by lemma 2.13 $cf(\aleph_\omega) = cf(\omega) = \omega < \aleph_\omega$ and hence by definition 2.14 $\aleph_\omega$ is singular. We set $\aleph_\alpha = \aleph_\omega$ and $\lambda = \aleph_0 < \aleph_\omega$ in Theorem 3.7 to conclude that there is a model $\mathcal{M}$ of $ZF$ (and hence of the weaker $Z$) in which $\hat{W}_{\aleph_0}$ holds, but on the other hand $AC_{cf(\omega)} \equiv AC_\omega \equiv AC_{\aleph_0}$ fails. Now, if a set $X$ is *Standard Infinite* then it is not finite and since $\hat{W}_{\aleph_0}$ holds, $card(X) \geq \aleph_0$. Thus, $SW_{\aleph_0}$ holds in $\mathcal{M}$. From Theorem 3.5 the statements $SW_{\aleph_0}$ and $SD$ are equivalent in $Z$, and thus in $\mathcal{M}$ the statement $SD$ holds, but the statement $AC_{\aleph_0}$ fails. That is what we wanted to prove. □

## **Appendix**

In this Appendix we prove auxiliary results mentioned in the paper. At first we provide a sketch of another proof of the fundamental result $SInf \Rightarrow CInf$ of Theorem 3.1. This proof is in $ZFC$ and it uses explicitly the *Axiom of Replacement* and the full version of the *Axiom of Choice* (see Ch. 8, Theorem 1.4 of [7]).

**Lemma A1**: ( $SInf \Rightarrow CInf$ ): *If $A$ is a standard infinite set, then $A$ contains a countable subset.*

Sketch of proof ($ZFC$): Let $A$ be a standard infinite set. From the Axiom of Choice ($AC$ AC), equivalently from the Well Ordering Theorem (see Ch.8 in [7]), $A$ is 'well orderable set' and hence using the *Axiom of Replacement* it is isomorphic with a unique infinite ordinal number $\Omega$ (see 6.3.1 in [7]). Since $\omega$ is the *least infinite ordinal*, $\Omega \geq \omega$. Now using the principle of transfinite recursion and the Axiom of Replacement (once more) we conclude that the set $A$ can be well ordered in a transfinite sequence of length $\Omega$, thus, $A = \langle a_\alpha : \alpha < \Omega \rangle$. From the defining properties of ordinals and well-ordered sets, $C = \langle a_\alpha : \alpha < \omega \rangle$ is an initial segment of $A$. Hence, the $Range(C) = \{a_\alpha : \alpha < \omega\}$ is a countable subset of $A$ (see also definition 1.2). □

In the next lemma we sketch a proof of *Pigeonhole Principle* that is used in proving $DInf \Rightarrow SInf$ in Theorem 3.1.

**Lemma A2**: **Pigeonhole Principle**. *If $A$ is a finite set, then there is no "one to one function" $f : A \to A$ onto a proper subset $B \subset A$. (Equivalently, if $A$ is a finite set, then every 'one to one' function $f : A \to A$ is also 'onto' $A$, and hence it is a bijection.)*

Sketch of Proof: We claim that to prove Pigeonhole Principle it is enough to show the following claim:

Claim: If $m \in \mathbb{N}$ and if $g : [0, m) \mapsto [0, m)$ is any 'one to one' function, then it is also an 'onto' function. (From the constructive definition of $\mathbb{N}$, $[0, m)$ is a subset of $\mathbb{N}$ that is identified with $m$).

Indeed, assuming the claim, consider any 'one to one' function $f : A \mapsto A$. Let $\pi : A \mapsto [0, m)$ be the bijection that witnesses the finiteness of $A$ for some $m \in \mathbb{N}$ (see definition 1.1). Then the function $g = \pi \circ f \circ \pi^{-1} : [0, m) \mapsto [0, m)$ is well defined and it is 'one to one'. Hence, form the claim it is also onto $[0, m)$. Thus, $f = \pi^{-1} \circ g \circ \pi$ is a bijection, as a composition of bijections, and hence it is onto $A$. Now the proof of the mentioned claim is a well-known and standard application of the Induction Principle, (see 5.27 of [10]). □